\documentclass[11pt]{amsart}
\usepackage{amsmath,amssymb}
\usepackage{setspace}
\usepackage{natbib}

\title{In Defence of the Material Conditional}
\author{Alexander V. Gheorghiu}
\address{University of Southampton, UK \and University College London, UK}
\email{a.v.gheorghiu@soton.ac.uk}

\begin{document}

\begin{abstract}
    The material conditional has long been charged with paradox. Defined truth-functionally, it renders true any conditional whose antecedent is false or consequent true --- hence, seemingly absurd statements such as `If unicorns exist, then  $2+2=4$'. This has been taken as proof that the connective cannot capture the meaning of ordinary if–then sentences, which appear to imply a causal or evidential link. I argue, by contrast, that the paradoxes arise from a confusion of what it expresses caused by cognitive biases. The material conditional properly belongs to the class of indicative, not subjunctive, conditionals --- those that register patterns of co-variation rather than counterfactual dependence. When understood as a formal device marking entailment under a background theory, it faithfully represents a mode of reasoning essential to science itself: correlation without causation. The faults ascribed to it, therefore, are not flaws in meaning or standard use, but as misapplication and misreading.
\end{abstract}


\maketitle
Few constructions in logic have provoked as much suspicion as the material conditional, $A \supset B$. This simple device, introduced by \citet{Frege1967}, has the virtue of formal elegance but the vice of indifference: in its truth-functional reading, it renders true many if–then statements that seem intuitively absurd. \citet{WhiteheadRussell1910} observed that it is equivalent to the statement `not-$A$ or $B$'. \citet{Pap1949} critiques this as an unpalatable result that `a true proposition is implied by any proposition, and a false proposition implies any proposition (the so-called paradoxes of material implication)' an outcome which, though tautologically correct, offends the ordinary sense in which if connects antecedent and consequent. On such grounds, \citet{Lewis1912} dismissed the connective as unfit to express genuine implication, proposing instead his modal notion of strict implication, under which ‘if $A$, then $B$’ holds only if $B$ follows of necessity from $A$.

Others have pressed the same objection from a pragmatic direction. \citet{Strawson1952} argued that an if-sentence conventionally suggests that its antecedent provides some reason for its consequent, and that, for this reason, we hesitate to call such statements true or false in abstraction. \citet{Ramsey1925}, in a similar spirit, remarked that if one already knows the antecedent to be false, such conditionals `cease to mean anything' --- the question whether `if $A$, then $B$?' simply does not arise. But the material conditional declares all such if-then true `by vacuity'. \citet{Grice1989} sought to defend it by relegating the divergence to pragmatics: it captures the literal truth-conditions, while our reluctance to assert conditionals with false antecedents arises from conversational implicature. \citet{Edgington1995} and others have noted, this defence sits uneasily with how we actually evaluate conditionals. Were material implication our ordinary standard we would be `intellectually disabled', unable to distinguish believable from unbelievable if-then statements.

I would like to re-examine these long-standing objections. Rather than rehearse them, I observe that they fail to understand what the material implication corresponds to in natural language. For despite its apparent paradoxes, the connective remains indispensable in formal reasoning. I shall argue that many of its supposed defects arise from a confusion between two notions of truth. Once that distinction is kept clear, the alleged paradoxes simply dissolve.
\medskip

First, I argue that the material conditional does capture a standard and legitimate use of if–then constructions. In both scientific and everyday reasoning, we often use conditionals to express patterns rather than causal relations. A familiar example comes from \citet{Hempel1965}:
\begin{quote}
If the barometer needle falls, then there will be a storm.
\end{quote}
No one imagines that the barometer causes the storm. The barometer falls because the air pressure drops; the same drop in pressure also produces the storm. The conditional thus expresses a systematic correlation under common conditions, not a causal chain.

Now suppose that on some particular day we are seated before a barometer and the needle is steady. The antecedent of the conditional is false, but does that fact bear on the truth of the statement? Of course does not: the conditional remains true because it expresses a general relation between pressure and weather, not a particular episode. 

Scientific instruments and diagnostic tests routinely rely on this very pattern: a variable $A$ co-occurs with $B$ because both are effects of a common cause $C$. When $A$ does not occur, the conditional $A \supset B$ remains true, since its purpose is to record the correlation, not the mechanism.   Against Ramsey, then, we may note that such conditionals do not `cease to mean anything' simply when the antecedent is false. In these cases, the implication is a test: we detect $A$ to infer $B$, because both arise in the same context. We do not expect our tests to be causes; to reject the material conditional because it fails to encode causation is therefore to overlook the very abstraction that makes scientific generalisation possible. Science students are taught that `correlation does not imply causation', perhaps textbooks on classical logic should add similar clarification.

In linguistic terms, these kind of if-then statements are \emph{indicative conditionals}. They express a standing or recurrent connection between antecedent and consequent. By contrast, \emph{subjunctive  conditionals} describe how things would stand under an imagined or  counterfactual situation --- for example, `If the barometer \emph{were} to fall, then there \emph{would} be a storm'. 

The distinction between indicative and subjunctive conditionals has, of course, been central since Stalnaker's work (\citeyear{Stalnaker1968,Stalnaker1975}) in which he formalised the contrast by introducing a possible-world semantics in which the truth of an indicative conditional depends on what holds at the closest contextually relevant world. \cite{Lewis1973} extended this framework to counterfactuals, giving the subjunctive reading its now-canonical semantics. On such approaches, every conditional receives a modal analysis, and the material conditional appears as a degenerate or impoverished fragment of a richer logic. Later accounts such as Jackson’s (\citeyear{Jackson1979,Jackson2006}) approached the indicative–subjunctive divide in terms of assertibility and contextual relevance, emphasising that our willingness to assert a conditional depends on the evidential relation between antecedent and consequent.

However, a more fine-grained linguistic analysis suggests a narrower correspondence. These scientific examples match a subset of indicative conditionals often described as \emph{lawlike} or \emph{generic} --- sentences that express regular correlations or generalisations over events. As \cite{Bennett2003} observes in his comprehensive taxonomy, lawlike indicatives constitute a distinct and well-established class of conditionals. A lawlike indicative conditional of the form `if $A$, then $B$' is typically paraphrasable as `$B$ whenever $A$'. My claim is that the material conditional abstracts these lawlike indicatives in the sense that its semantics captures their meaning. Here I have only shown that such conditionals are pervasive in ordinary language; below I show how the material conditional abstracts them. Before turning to the arguments, I will first consider why certain constructions appear absurd. 
\medskip

Second, I argue that the apparent absurdities of the material conditional arise from misreading  lawlike indicative conditionals as saying something other than intended. These misreadings are often justified by attempting to slot statements whose context is unclear into usual modes of discourse; but they are misreadings nonetheless. As a representative example, consider the apparently nonsensical statement:
\begin{quote}
If unicorns exist, then $2+2=4$.
\end{quote}
This example is pathological for several reasons.
\begin{enumerate}
\item unicorns do not, in fact, exist;
\item the consequent, $2+2=4$, is unconditionally true;
\item there is no discernible connection between the existence of unicorns and the truth of an arithmetic identity.
\end{enumerate}
The sense of absurdity arises really from (3): the clauses are drawn from unrelated domains, so the conditional conveys no intelligible pattern or relation. Once that lack of connection is recognised, the impression that the statement is `true by vacuity' under (1) and (2) is revealed as a by-product of misunderstanding. When the conditional is interpreted in its proper lawlike indicative sense it becomes quite clear. 

On (1), we can recognise the force of Ramsey’s familiar objection: that indicative conditionals are meaningless when their antecedent is known to be false. As noted earlier, this objection misfires. The difficulty here is not that unicorns do not exist, but that their existence is not something we can meaningfully vary --- it is not, as it were, a condition that can be `turned on and off'. The antecedent thus fails to function as a test variable, and the conditional consequently strikes us as empty. Faced with this, we instinctively shift to a subjunctive reading: if unicorns were to exist, then it would be that $2+2=4$. Interpreted that way, the sentence acquires a counterfactual or even ironic flavour --- roughly, yes, $2+2=4$, when pigs fly. This shows that the difficulty lies not with the statement itself but with our manner of reading it. Stripped of that pragmatic bias it is entirely coherent. Indeed, everyone would agree that `$2+2=4$ whenever unicorns exist', for $2+2=4$ holds under all conditions whatsoever.

On (2), we encounter Pap’s objection: that a true proposition is implied by any proposition, a result which, as he put it, `offends' our ordinary understanding of if–then constructions. The source of discomfort here again lies not in the truth of $2+2=4$, but in reading the conditional as suggesting that we imagine circumstances under which $2+2=4$ might fail. Since we cannot coherently do so, the statement seems nonsensical.

This reaction, however, again stems from a misplaced subjunctive reading. It seems to say that the fact that $2+2=4$ depends on whether or not unicorns exist --- plainly absurd. To see this, compare the statement `If Peano Arithmetic is true, then $2+2=4$'. No one finds that problematic, though the consequent is unchanged, because the antecedent and consequent belong to the same explanatory field. But read properly as a lawlike indicative conditional, the claim is much weaker: it only says `$2+2=4$ whenever unicorns exist' and that is quite reasonable for $2+2=4$ is always true. The conditional is quite sensible in this reading, if somewhat gratuitous.

On (3), we turn to Strawson’s observation that if–then constructions ordinarily suggest some intelligible link between antecedent and consequent. As noted, this link need not be causal or temporal. But then what link is actually expressed by lawlike indicative conditionals? In this case, it is co-variation under a background theory. Within the framework of meteorology, for example, we say that a storm occurs whenever the barometer needle falls for they have the common cause of low pressure; within a mythological framework, one might equally say that a storm occurs whenever Zeus is angry. Each conditional expresses a lawlike relation relative to its governing theory. It is in making this analysis precise that we require the semantics of the material conditional. \medskip

I claim that the material conditional abstracts lawlike indicative conditionals. By this I mean that the semantics of the material conditional is an accurate mathematical rendition of the intended meaning of lawlike indicative conditionals. The essential observation is that in ordinary language we often oscillate between two different notions of `truth' but in logic we distinguish them. That some ordinary language statements appear absurd when read as material conditionals comes down to reading them with the wrong notion of truth. 

The first notion of truth is \emph{truth within a theory}.  When we say that a conditional $A \supset B$ is true relative to some background assumptions~$\Gamma$, we mean that it can be derived from them --- that the judgement 
\[
\Gamma \vdash A \supset B
\]
is warranted. Here `true' does not mark correspondence with some particular fact but conformity with a network of assumptions --- the way a statement fits within a framework of reasoning. 

Logicians rarely use the word `true' for this sort of relation; they prefer to say that $A \supset B$ is \emph{provable} or \emph{derivable} from~$\Gamma$. But in ordinary and scientific language we do sometimes speak this way.  
We say, for example, that a claim is `true of Newtonian mechanics' or `true in the theory of evolution'. When pushed, we might clarify that by true when mean `it follows from' or `it is a consequence from', but otherwise we often employ the predicate as stated it.  

The second notion of truth is \emph{truth within a model} as given by \cite{Tarski1933}.  This is, perhaps, the more familiar notion in logic based on interpretation.   A statement is true in a model~$\mathfrak{M}$, written $\mathfrak{M} \models A \supset B$, when the model represents the antecedent as false or the consequent as true:
\[
\mathfrak{M} \models A \supset B 
\quad\text{iff}\quad
\mathfrak{M} \not\models A \text{ or } \mathfrak{M} \models B.
\]
 This notion of truth is the logician’s counterpart to the ordinary use in which we say that something is `true of arithmetic' referring to, say, the standard model $\mathfrak{A}$ of Peano Arithmetic (PA).  The semantic conception thus gives formal expression to a very familiar idea that a statement is true when what it says is borne out by the system or world it describes.

The two notions of truth are related through  Tarski's (\citeyear{Tarski1936}) account of logical consequence:
\[ 
\Gamma \vdash A \supset B \qquad \mbox{iff} \qquad \mbox{for any model $\mathfrak{M}$, if $\mathfrak{M} \models C$ for $C \in \Gamma$, then $\mathfrak{M} \models A \supset B$} 
\]
The distinction, however, is easier to blur than to maintain.  In logic, we formally delineate these notions of truth; but in general parlance, we rarely stop to clarify which notion of truth we have in mind. 

Of course, we must also be careful in our choice of background theory.  One could, in principle, erase the distinction between theory and model by constructing a theory that simply asserts everything true in a given model:
\[
\Omega(\mathfrak{M}) := 
\{\,C \mid \mathfrak{M} \models C\,\} \cup 
\{\,\neg C \mid \mathfrak{M} \not\models C\,\}.
\]
Such a theory would reproduce the model point for point, but it would not be a theory in the relevant sense.  No serious theory of meteorology, for instance, contains the weather conditions of a particular day within its axioms.  

When we consider the status of lawlike indicative conditionals, we do so not by the status of the antecedent or consequence in a particular circumstance but whether or not it holds as a general law. Mathematically, this is expressed as truth-in-a-theory. Hempel’s barometer conditional, for example, is true because it holds within our theory of meteorology. Whether, on some particular day, the barometer needle happens to fall is beside the point. Such local details belong to a `model' --- in the logician's sense, not the scientist's --- that specifies the actual weather and the needle’s position. They do not affect the `truth' of the conditional itself, which is secured by the theory linking air pressure to storms. The purpose of truth-in-a-model here is not to determine the truth the implication statement but to explicate the  co-variation it expresses: ranging over models of the theory, the conclusion is true whenever the antecedent is true. This is the co-variation expressed by lawlike indicative conditionals, which is independent from their standing as true or false.

How has this confusion between notions of truth historically entered our analysis? Typically, when we introduce the material implication  we do so through its truth-functional behaviour using, say, truth-tables. We then gives ordinary language examples such as those above as examples of such implication. The problem is that what it means to evaluate such sentences as true or false is not the same as what makes the material conditional true or false in a given model.  It is this modulation of the notion of truth between theory and example that makes the material conditional reading of certain if-then statements seem absurd.

Let us then return to the apparent nonsense of `if unicorns exist, then $2+2=4$' and see how this understanding resolves the seeming absurdity. When conditionals are presented without an appropriate theory, we easily slip from truth-in-a-theory to truth-in-a-model to make sense of it. Such slippages are understandable: they reflect the ease with which our linguistic habits drift between different conceptions of truth.  But they are nonetheless mistakes, for they alter the meaning of the statement and change what we mean in calling such statements true or false.   The difficulties noted in (1) and (2) stem precisely from this slippage. No reasonable background theory seemingly contains both the metaphysical status of unicorns and the arithmetic of integers;  consequently, we evaluate the sentence by reference to the actual facts of today. 

It is the missing context which is supposed to link the existence of unicorns and arithmetic to answer (3). In the absence of a context,  we cannot say that the statement `if unicorns exist, then $2+2=4$' is either true or false. To do so presupposes a background theory relative to which the evaluation is made. Without a context, calling the statement true or false is not so much wrong as specious. Therefore, it is quite incorrect to suggest that reading such pathological if-then statements as material implications renders them true. Accordingly, the apparent absurdities of material implication vanish. 

We might attempt to read such statements relative to some context such as PA. This is dubious for the existence of unicorns lies beyond the scope of the language of arithmetic in which PA is formulated. Nonetheless, properly contextualised in this way, we don't need to ask `if unicorns exist, then $2+2=4$’ simpliciter is true, but rather if the statement
\[
\text{It is true of PA that $2+2=4$ whenever unicorns exist.}
\]
is true. There is nothing absurd about this, it is merely gratuitous for it simply observes that in any model of PA in which unicorns exist is also a model of $2+2=4$. 
\medskip

The material conditional has survived its critics not through stubborn tradition but because it captures standard uses of if-then constructions in ordinary language.   It does not capture every such construction --- natural language is far too various for that—but it captures precisely those that matter most in mathematics and science: the lawlike indicative conditionals.  These express not a causal link between antecedent and consequent, but a co-variation of truth across admissible situations.  

The confusion it has caused has, perhaps, two sources.  First, logic textbooks often begin with a discussion of truth rather than of consequence, encouraging the reader to think in terms of truth-values rather than inferential roles.  Second, when confronted with conditionals in natural language, we naturally drift toward whatever interpretation feels most familiar --- reading indicatives as subjunctives, or appealing to the wrong notion of truth --- even when those readings are, strictly speaking, mistaken.  Clarifying this distinction restores to the material conditional its proper status: not as a defective surrogate for everyday reasoning, but as the precise formal expression of a lawlike relation between antecedent and consequent. 

\bibliographystyle{apalikefull}
\bibliography{bib}
\end{document}